\documentclass[12pt]{article}

\usepackage{theorem,amssymb}

\advance \topmargin by -\headheight
\advance \topmargin by -\headsep
\evensidemargin \oddsidemargin
\author{J.-P. Allouche \\
CNRS, LRI, B\^atiment 490 \\
F-91405 Orsay Cedex (France) \\
{\tt allouche@lri.fr} \\
}

\title{On a conjecture of Deutsch, Sagan, and Wilson}

\date{ }

\def \proof{\bigbreak\noindent{\it Proof.\ \ }}

\def \endpf{{\ \ $\Box$ \medbreak}}

\newtheorem{proposition}{Proposition}
\theorembodyfont{\rm}

\newtheorem{remark}{Remark}

\begin{document}

\maketitle

\begin{abstract}
We prove a recent conjecture due to Deutsch, Sagan, and Wilson stating
that the finite sequence obtained from the first $p$ central trinomial 
coefficients modulo $p$ by replacing nonzero terms by $1$'s is palindromic, 
for any prime number $p \geq 5$.
\end{abstract}

\section{Introduction}
In the recent paper \cite{DeuSag} Deutsch and Sagan study several combinatorial 
sequences reduced modulo prime numbers. They are in particular interested in 
the values modulo $p$ of the central trinomial coefficients. Let us recall that
the $n$th central trinomial coefficient is defined as the largest coefficient 
in the expansion of the polynomial $(1 + x + x^2)^n$. Deutsch and Sagan make 
the following conjecture \cite[Conjecture 5.8]{DeuSag} (also stated by Wilson,
see \cite{DeuSag}): for each prime $p \geq 5$ and for each $j < p$, the 
number $T_j$ is divisible by $p$ if and only if $T_{p-1-j}$ is divisible by $p$.
We give here an elementary proof of this conjecture.

\section{The central trinomial coefficients modulo a prime}

In this section we first recall a classical result 
(see \cite[sequence A002426]{Sloane} for example).

\begin{proposition}\label{classical}
The generating function $F$ of the central trinomial coeficient satisfies:
$$
F(x) := \sum_{n \geq 0} T_n x^n = \frac{1}{\sqrt{1 - 2x -3x^2}}\cdot
$$
\end{proposition}

\begin{remark}\label{Lucas}
Note that, using \cite[Theorem 6.4]{Allo-trans} (see also \cite{AGS}), this 
implies that the sequence $(T_n)_{n \geq 0}$ satisfies: if $p$ is prime and 
if the base $p$ expansion of $n$ is $n = \sum n_j p^j$, then 
$T_n \equiv \prod T_{n_j} \bmod p$, which is Theorem 4.7 of \cite{DeuSag}.
\end{remark}

An easy consequence of Proposition~\ref{classical} is the following statement.

\begin{proposition}\label{polynomial}
Let $p$ be an odd prime. Then we have the following identity
$$
\sum_{0 \leq n \leq p-1} T_n x^n \equiv (1 - 2x -3x^2)^{\frac{p-1}{2}} \bmod p.
$$
\end{proposition}

\proof From Proposition~\ref{classical} we have
\begin{equation}\label{eq1}
(1 - 2x - 3x^2)^{\frac{p-1}{2}} \left(\sum_{n \geq 0} T_n x^n\right)^{p-1} 
\equiv 1 \bmod p.
\end{equation}
On the other hand, using the ``$p$-Lucas property'' recalled in Remark~\ref{Lucas}
above, we have
$$
\begin{array}{lll}
\displaystyle\sum_{n \geq 0} T_n x^n &\equiv 
\displaystyle\sum_{0 \leq j \leq p-1} \sum_{n \geq 0} T_{pn+j} x^{pn+j}
&\equiv \displaystyle\sum_{0 \leq j \leq p-1} \sum_{n \geq 0} T_nT_j x^{pn} x^j \\
&\equiv \displaystyle\sum_{0 \leq j \leq p-1} T_j x^j \sum_{n \geq 0} T_n x^{pn}
&\equiv \left(\displaystyle\sum_{0 \leq j \leq p-1} T_j x^j\right)
\displaystyle\left(\sum_{n \geq 0} T_n x^n\right)^p \bmod p
\end{array}
$$
which yields
\begin{equation}\label{eq2}
\left(\sum_{0 \leq j \leq p-1} T_j x^j\right)
\left(\sum_{n \geq 0} T_n x^n\right)^{p-1}
\equiv 1 \bmod p.
\end{equation}
Comparing Equations \ref{eq1} and \ref{eq2} finishes the proof. \endpf

\section{Proof of the conjecture}

We first prove a proposition on the nonzero coefficients of a quadratic 
polynomial raised to an integer power.

\begin{proposition}\label{powers}
Let $1 + ax + bx^2$ be a polynomial with coefficients in a commutative 
field $K$, with $b \neq 0$. Let $k$ be a positive integer. Then, noting
$(1 + ax + bx^2)^k := \sum_{0 \leq j \leq 2k} \alpha_j x^j$, we have
$\alpha_j = 0$ if and only if $\alpha_{2k-j} = 0$.
\end{proposition}

\proof We write 
$$
\sum_{0 \leq j \leq 2k} \alpha_j b^{k-j} x^{2k-j}
= b^k x^{2k} \displaystyle\sum_{0 \leq j \leq 2k} \alpha_j \left(\frac{1}{bx}\right)^j 
= b^k x^{2k} \left(1 + \displaystyle\frac{a}{bx} + \frac{b}{(bx)^2}\right)^k
= (bx^2 + ax + 1)^k.
$$
But the sum on the left can also be written 
$\sum_{0 \leq j \leq 2k} \alpha_{2k-j}  b^{j-k} x^j$;
thus, for all $j \in [0, 2k]$, we have $\alpha_{2k-j} = b^{k-j} \alpha_j$
which implies our claim. \endpf

\bigskip

As an immediate corollary, we get a proof of the conjecture of Deutsch, Sagan, 
and Wilson.

\bigskip

\noindent
{\bf Theorem} {\it For any prime $p \geq 5$, for any $j \in [0, p-1]$, the sequence 
of central trinomial coefficients $(T_n)_{n \geq 0}$ satisfies
$$
p \ | \  T_j \ \mbox{\rm if and only if} \ p \ | \ T_{p-1-j}.
$$
}

\medskip

\proof Apply Proposition~\ref{powers} with $K := {\mathbb Z}/p {\mathbb Z}$
(the finite field with $p$ elements) and $1+ax+bx^2 := 1-2x-3x^2$, and use 
Proposition~\ref{polynomial}. \endpf

\begin{remark}
The reader can check that the proof of the Theorem above readily generalizes to
proving the following. (Hint: use \cite[Theorem 2 and its proof]{AGS}.)

\noindent
{\it Let $(R_n)_{n \geq 0}$ be a sequence of integers, such that there exists
a polynomial of degree $2$ with integer coefficients $P(x) := 1 + ax + bx^2$
such that $\sum_{n \geq 0} R_n x^n = (P(x))^{-1/2}$. Then, for all primes $p$
such that $p$ does not divide $3R_1^2 - 2R_2$ and for all $j \in [0, p-1]$, we 
have 
$$
p \ | \ R_j \ \mbox{\rm if and only if} \ p \ | R_{p-1-j}.
$$
}

In particular if $(R_n)_{n \geq 0}$ is the sequence of central Delannoy numbers
(see \cite[sequence A001850]{Sloane}), then for all primes $p$ and for 
all $j \in [0, p-1]$, we have
$$
p \ | \ R_j \ \mbox{\rm if and only if} \ p \ | R_{p-1-j}.
$$
Note that the $p$-Lucas property for this sequence is a consequence of 
 \cite[Theorem 6.4]{Allo-trans} (see also \cite{AGS}) and of the fact that
the generating function for the central Delannoy numbers is equal to
$(1 - 6x + x^2)^{-1/2}$ (see \cite[sequence A001850]{Sloane} for example);
it is also proven in \cite{DeuSag} and in \cite{ELY}. A nice paper on sequences 
having the $p$-Lucas property is \cite{McIntosh}.
\end{remark}

\noindent
{\bf Addendum:} the result was proved before almost in the same way by 
Tony D. Noe: On the Divisibility of Generalized Central Trinomial Coefficients,
Journal of Integer Sequences, Vol. 9 (2006), Article 06.2.7
          
\noindent
{\tt http://www.cs.uwaterloo.ca/journals/JIS/VOL9/Noe/noe35.html}

\end{document}